\documentclass[12pt]{amsart}
\def\NZQ{\Bbb}               

\def\QQ{{\NZQ Q}}
\def\ZZ{{\NZQ Z}}

\def\frk{\frak}               

\def\Phi{{\frk n}}
\def\Phi{{\frk N}}
\def\opn#1#2{\def#1{\operatorname{#2}}} 
\opn\chara{char}
\opn\length{\ell}
\opn\pd{pd}
\opn\rk{rk}
\opn\projdim{proj\,dim}
\opn\injdim{inj\,dim}
\opn\rank{rank}
\opn\depth{depth}
\opn\grade{grade}
\opn\height{height}
\opn\embdim{emb\,dim}
\opn\codim{codim}

\opn\Tr{Tr}
\opn\bigrank{big\,rank}
\opn\superheight{superheight}\opn\lcm{lcm}
\opn\trdeg{tr\,deg}%
\opn\reg{reg}
\opn\lreg{lreg}
\opn\ini{in}
\opn\lpd{lpd}
\opn\div{div}
\opn\Div{Div}
\opn\cl{cl}
\opn\Cl{Cl}
\opn\Spec{Spec}
\opn\Supp{Supp}
\opn\supp{supp}
\opn\Sing{Sing}
\opn\Ass{Ass}
\opn\Ann{Ann}
\opn\Rad{Rad}
\opn\Soc{Soc}
\opn\Im{Im}
\opn\Ker{Ker}
\opn\Coker{Coker}
\opn\Am{Am}
\opn\Hom{Hom}
\opn\Tor{Tor}
\opn\Ext{Ext}
\opn\End{End}
\opn\Aut{Aut}
\opn\id{id}

\opn\nat{nat}
\opn\pff{pf}
\opn\Pf{Pf}
\opn\GL{GL}
\opn\SL{SL}
\opn\mod{mod}
\opn\ord{ord}
\opn\Gin{Gin}
\opn\aff{aff}
\opn\con{conv}
\opn\relint{relint}
\opn\st{st}
\opn\lk{lk}
\opn\cn{cn}
\opn\core{core}
\opn\vol{vol}
\opn\link{link}
\opn\star{star}
\opn\com{com}
\opn\gr{gr}

\def\pot#1#2{#1[\kern-0.28ex[#2]\kern-0.28ex]}

\opn\dirlim{\underrightarrow{\lim}}
\opn\inivlim{\underleftarrow{\lim}}
\let\to=\rightarrow

\def\Implies{\ifmmode\Longrightarrow \else
       \unskip ${}\Longrightarrow{}$\ignorespaces\fi}
\def\implies{\ifmmode\Rightarrow \else
       \unskip ${}\Rightarrow{}$\ignorespaces\fi}
\def\iff{\ifmmode\Longleftrightarrow \else
       \unskip ${}\Longleftrightarrow{}$\ignorespaces\fi}

\let\:=\colon
\newtheorem{Theorem}{Theorem}[section]
\newtheorem{Lemma}[Theorem]{Lemma}
\newtheorem{Corollary}[Theorem]{Corollary}
\newtheorem{Proposition}[Theorem]{Proposition}
\newtheorem{Remark}[Theorem]{Remark}

\newtheorem{Example}[Theorem]{Example}

\newtheorem{Definition}[Theorem]{Definition}

\let\epsilon\varepsilon
\let\phi=\varphi
\let\kappa=\varkappa
\def\qed{\ifhmode\textqed\fi
     \ifmmode\ifinner\quad\qedsymbol\else\dispqed\fi\fi}
\def\textqed{\unskip\nobreak\penalty50
      \hskip2em\hbox{}\nobreak\hfil\qedsymbol
      \parfillskip=0pt \finalhyphendemerits=0}
\def\dispqed{\rlap{\qquad\qedsymbol}}

\opn\dis{dis}
\def\pnt{{\raise0.5mm\hbox{\large\bf.}}}

\opn\lex{lex}
\opn\rev{rev}
\opn\Lex{Lex}
\opn\GL{GL}
\opn\initial{in}
\opn\gin{gin}

\newcommand{\Mac}{ ^{\mathrm{MG}(n-1)}}
\newcommand{\Macd}{ ^{\mathrm{MG}(d)}}
\newcommand{\sat}{{\mathrm{sat}}}

\begin{document}
\title{Gotzmann ideals of the polynomial ring}
\author{Satoshi Murai and Takayuki Hibi}
\address{Satoshi Murai, Department of Pure and Applied Mathematics,
Graduate School of Information Science and Technology,
Osaka University, Toyonaka, Osaka 560-0043, Japan}
\email{s-murai@ist.osaka-u.ac.jp}
\address{Takayuki Hibi, Department of Pure and Applied Mathematics,
Graduate School of Information Science and Technology,
Osaka University, Toyonaka, Osaka 560-0043, Japan}
\email{hibi@math.sci.osaka-u.ac.jp}
\date{}
\thanks{The first author is supported by JSPS Research Fellowships for Young Scientists.}
\maketitle
\begin{abstract}
Let
$A = K[x_1, \ldots, x_n]$
denote the polynomial ring in $n$ variables
over a field $K$.
We will classify all the Gotzmann ideals of $A$ with
at most $n$ generators.
In addition,
we will study Hilbert functions $H$ for which all homogeneous ideals of $A$ with the Hilbert function $H$
have the same graded Betti numbers. These Hilbert functions will be called inflexible Hilbert functions.
We introduce the notion of segmentwise critical Hilbert function and show that
segmentwise critical Hilbert functions are inflexible.
\end{abstract}

\section{Introduction}
Let
$K$ be an arbitrary field and
let $A = K[x_1, \ldots, x_n]$
denote the polynomial ring in $n$ variables
over $K$
with each $\deg x_i = 1$.
For a homogeneous ideal $I$ of $A$,
the \textit{Hilbert function $H(I,-): \mathbb{Z}_{\geq 0} \to \mathbb{Z}_{\geq 0}$ of $I$} is the numerical function
defined by $H(I,t)=\dim_K I_t$, where $I_t$ is the homogeneous component of degree $t$ of $I$.
Work with the lexicographic order $<_{\lex}$
on $A$ induced by the ordering
$x_1 > x_2 > \cdots > x_n$.
Recall that a set $V$ of monomials in $A$ is said to be
{\em lexsegment} if,
for monomials $u$ and $v$ of $A$
with $u \in V$, $\deg u = \deg v$ and
$u <_{\lex} v$, one has $v \in V$.
Also, a monomial ideal $I$ of $A$ is called a \textit{lexsegment ideal} if 
the set of monomials in $I$ is lexsegment.
Let $I$ be a homogeneous ideal of $A$
and $I^{\lex}$ the (unique) lexsegment ideal
(\cite{Big} and \cite{Hul})
with the same Hilbert function as $I$.
A {\em Gotzmann ideal} introduced in \cite{HerHi} is
a homogeneous ideal $I$ for which
the number of minimal generators of $I$
is equal to that of $I^{\lex}$.
In this paper,
we classify all the Gotzmann
ideals of $A$ generated by at most $n$ homogeneous
polynomials.

A numerical function $H(-) : \mathbb
{Z}_{\geq 0} \to \mathbb{Z}_{\geq 0}$ is said to be {\em critical} \cite{MuHi} if it is equal to the Hilbert function of
a lexsegment ideal of $A$ which has at most $n$ generators.
Let $1 \leq s \leq n$ and
$f_1, \ldots, f_s$ homogeneous polynomials
with
\[
f_i \in K[x_i, x_{i+1}, \ldots, x_n]
\]
for each $1 \leq i \leq s$
and with $\deg f_s > 0$.
In \cite{Got}
the ideal
$I_{(f_1, \ldots, f_s)}$ of $A$
defined by
\begin{eqnarray}
\label{canonicalcritical}
I_{(f_1, \ldots, f_s)}
= (f_1x_1, f_1f_2x_2, \ldots,
f_1 f_2 \cdots f_{s-1} x_{s-1},
f_1 f_2 \cdots f_s)
\end{eqnarray}
was introduced.
A homogeneous ideal $I$ of $A$
is called {\em canonical critical}
if $I = I_{(f_1, \ldots, f_s)}$
for some homogeneous polynomials
$f_1, \ldots, f_s$
with
$f_i \in K[x_i, x_{i+1}, \ldots, x_n]$
for each $1 \leq i \leq s$
and with $\deg f_s > 0$,
where $1 \leq s \leq n$.
Our first result is

\begin{Theorem}
\label{gotzmann}
Given a homogeneous ideal $I$ of
$A = K[x_1, \ldots, x_n]$,
the following conditions are equivalent:
\begin{enumerate}
\item[(i)]
$I$ has a critical Hilbert function;
\item[(ii)] there exists
a linear transformation
$\varphi$ on $A$
such that $\varphi(I)$ is a canonical
critical ideal;
\item[(iii)] $I$ is a Gotzmann ideal
generated by at most $n$ homogeneous polynomials.
\end{enumerate}
\end{Theorem}

It was shown in \cite[Corollary 1.4]{HerHi} that
a homogeneous ideal $I$ is Gotzmann if and only if
$I$ and $I^{\lex}$ have the same graded Betti numbers.
Thus Theorem \ref{gotzmann} shows that if a homogeneous ideal
has a critical Hilbert function $H$,
then the graded Betti numbers
are determined by $H$.
In the second part of this paper,
we generalize this fact.
We introduce two classes of Hilbert functions
$H$
for which all ideals with the Hilbert function $H$
have the same graded Betti numbers (Theorems \ref{2-5} and \ref{Peeva}).

We establish fundamental properties of critical
functions and canonical critical ideals in Section $2$, and give
our proof of Theorem \ref{gotzmann}
in Section $3$.
Next, in Section $4$, we introduce some Hilbert functions $H$ such that
all ideals with the Hilbert function $H$ have the same graded Betti numbers.
Finally, in the appendix,
by using Theorem \ref{gotzmann},
we give a simple and purely algebraic proof
of \cite[Proposition 2]{Got}.

\section{Critical ideals}
In this section,
we establish fundamental properties of critical ideals.
First, we recall basic facts on Gotzmann ideals and lexsegment ideals.
Let $K$ be a field and
$A = K[x_1, \ldots, x_n]$
the polynomial ring in $n$ variables
over a field $K$
with each $\deg x_i = 1$.
Let $d\geq 0$ be an integer and $V$ a subspace of the $K$-vector space $A_d$.
We write $\Lex(V) \subset A_d$ for the $K$-vector space
spanned by the lexsegment set $L \subset A_d$ of monomials
with $|L|= \dim_K V$.
A famous result of Macaulay (see e.g., \cite[Corollary C4]{IK})
guarantees
\begin{eqnarray}
\label{MacIneqSpace}
 \dim_K (A_1 \cdot V) \geq \dim_K (A_1 \cdot \Lex(V)),
\end{eqnarray}
where $A_1 \cdot V = \{ \ell g: \ell \in A_1 \mbox{ and } g \in V\}$.
A $K$-vector space $V$ is called a \textit{Gotzmann space} if 
$\dim_K (A_1 \cdot V) = \dim_K (A_1 \cdot \Lex(V))$.
A homogeneous ideal $I$ of $A$ is said to be \textit{Gotzmann} if
$I_k$ is Gotzmann for all $k \geq 0$.
For any homogeneous ideal $I \subset A$,
let $I^{\lex} = \bigoplus_{k \geq 0} \Lex(I_k)$.
It follows from (\ref{MacIneqSpace}) that
$I^{\lex}$ is indeed a lexsegment ideal
having the same 
Hilbert function as $I$.
Clearly a homogeneous ideal $I$ of $A$ 
is Gotzmann if and only if
the number of minimal generators of $I$ is equal to that of $I^{\lex}$.

A monomial ideal $I$ of $A$
is said to be \textit{universal lexsegment} \cite{BNT}
if, for any integer $m \geq 0$, the ideal of $K[x_1,\dots,x_{n+m}]$
having the same generators as $I$ is a lexsegment ideal of $K[x_1,\dots,x_{n+m}]$.
It was shown in \cite[Corollary 1.3]{MuHi} that
a monomial ideal $I$ of $A$ is universal lexsegment
if and only if
$I$ is lexsegment and $|G(I)| \leq n$, where
$G(I)$ is the set of the minimal monomial
generators of $I$.
Moreover, it is known \cite[Proposition 1.2]{MuHi} that
a monomial ideal $I$ of $A$ is universal lexsegment
if and only if, for some $1 \leq s \leq n$
and for some nonnegative integers
$b_1, b_2, \ldots, b_s$,
one has
\begin{eqnarray}
\label{unilex}
I = (x_1^{b_1+1}, x_1^{b_1}x_2^{b_2+1}, \ldots,
x_1^{b_1} \cdots x_{s-1}^{b_{s-1}} x_s^{b_s+1}).
\end{eqnarray}
It follows from \cite[Proposition 1.5]{MuHi}
that the Hilbert function $H(I,t)$ of
the universal lexsegment ideal
$(\ref{unilex})$
is given by
\begin{eqnarray}
\label{hilbertfunction}
H(I,t) = {t - a_1 + n - 1 \choose n - 1}
+ \cdots +
{t - a_s + n - s \choose n - s},
\end{eqnarray}
where the sequence
$(a_1, a_2, \ldots, a_s)$ with
\[
0 < a_1 \leq a_2 \leq \cdots \leq a_s
\]
is defined by setting
\[
a_i = \deg x_1^{b_1} \cdots x_{i-1}^{b_{i-1}} x_i^{b_i+1},
\, \, \, \, \, \, \, \, \, \,
1 \leq i \leq s
\]
or equivalently,
\[
b_i = a_i - a_{i-1},
\, \, \, \, \, \, \, \, \, \,
1 \leq i \leq s
\]
with $a_0 = 1$.

Since a lexsegment ideal with a given Hilbert function is
uniquely determined, it follows that
a numerical function $H: \ZZ_{\geq 0} \to \ZZ_{\geq 0}$ is critical if and only if
there exists a sequence $(a_1, \ldots, a_s)$ of
integers with
$0 < a_1 \leq a_2 \leq \cdots \leq a_s$,
where $1 \leq s \leq n$, such that
$H(t)$ is written in the form
$(\ref{hilbertfunction})$.
A critical function of the form
$(\ref{hilbertfunction})$
will be called a critical function of
{\em type $(a_1, a_2, \ldots, a_s)$}.
For convenience, we say that a homogeneous ideal of $A$
is a \textit{critical ideal of type $(a_1, a_2, \ldots, a_s)$}
if its Hilbert function is the critical function
of type $(a_1, a_2, \ldots, a_s)$.

\begin{Lemma}
\label{Boston}
Let $1 < s \leq n$.
Fix homogeneous polynomials
$f_1, \ldots, f_{s-1}$ with each
$f_i \in K[x_i, \ldots, x_n]$.
Let $g \in K[x_s, \ldots, x_n]$
be a homogeneous polynomial
with $\deg g > 0$.
Then
\begin{eqnarray}
\label{notin}
f_1 f_2 \cdots f_{s-1} g
\not\in
(f_1x_1, f_1f_2x_2, \ldots,
f_1 f_2 \cdots f_{s-1} x_{s-1}).
\end{eqnarray}
\end{Lemma}

\begin{proof}
One has $(\ref{notin})$
if and only if
\begin{eqnarray}
\label{notinin}
f_2 \cdots f_{s-1} g
\not\in
(x_1, f_2x_2, \ldots,
f_2 \cdots f_{s-1} x_{s-1}).
\end{eqnarray}
Since
$f_2 \cdots f_{s-1} g
\in K[x_2, \ldots, x_n]$,
it follows that $(\ref{notinin})$
is equivalent to saying
\begin{eqnarray}
\label{notininin}
f_2 \cdots f_{s-1} g
\not\in
(f_2x_2, \ldots,
f_2 \cdots f_{s-1} x_{s-1}).
\end{eqnarray}
Now, working with induction on $s$,
the desired result $(\ref{notininin})$
is guaranteed.
\end{proof}

\begin{Lemma}
\label{directsum}
As a vector space over $K$
the ideal
$(\ref{canonicalcritical})$
is the direct sum
\begin{eqnarray}
\label{vectorspace}
I_{(f_1, \ldots, f_s)}
= \left(\bigoplus_{j=1}^{s-1}
(f_1 \cdots f_j x_j)K[x_j, \ldots, x_n] \right)
\bigoplus
(f_1 \cdots f_s)K[x_s, \ldots, x_n].
\end{eqnarray}
\end{Lemma}

\begin{proof}
Let $J = (f_1x_1, f_1f_2x_2, \ldots,
f_1 f_2 \cdots f_{s-1} x_{s-1})$.
Let $g = f_{s-1}x_{s-1}$.
Then
\[
J = (f_1x_1, \ldots,
f_1 f_2 \cdots f_{s-2} x_{s-2},
f_1 f_2 \cdots f_{s-2} g).
\]
Now, working with induction on $s$ yields
\[
J = \left(\bigoplus_{j=1}^{s-2}
(f_1 \cdots f_j x_j)K[x_j, \ldots, x_n] \right)
\bigoplus
(f_1 f_2 \cdots f_{s-2} g)K[x_{s-1}, \ldots, x_n].
\]
In other words, one has
\[
J = \bigoplus_{j=1}^{s-1}
(f_1 \cdots f_j x_j)K[x_j, \ldots, x_n].
\]
On the other hand, by using Lemma \ref{Boston},
it follows that
\[
(f_1 \cdots f_{s-1}f_s)K[x_s, \ldots, x_n]
\bigcap J = \{0\}.
\]
Clearly, $f_1 \cdots f_{s-1}f_s x_j \in J$
for each $j < s$.  Hence
$(\ref{vectorspace})$
follows.
\end{proof}

\begin{Corollary}
\label{Sydney}
Let $I_{(f_1, \ldots, f_s)}$
denote the ideal
$(\ref{canonicalcritical})$.
\begin{itemize}
\item[(a)]
$I_{(f_1, \ldots, f_s)}$
is a critical ideal of
type $(a_1, \ldots, a_s)$,
where
$a_i = \deg f_1 f_2 \cdots f_i x_i$,
$i = 1, \ldots, s - 1$, and
$a_s = \deg f_1 f_2 \cdots f_s$.
\item[(b)]
$I_{(f_1, \ldots, f_s)}$
is minimally generated by
\begin{eqnarray}
\label{minimally}
\{ f_1x_1, \ldots,
f_1 f_2 \cdots f_{s-1} x_{s-1},
f_1 f_2 \cdots f_s \}.
\end{eqnarray}

\item[(c)]
$I_{(f_1, \ldots, f_s)}$
is Gotzmann.
\end{itemize}
\end{Corollary}

\begin{proof}
The direct sum decomposition
$(\ref{vectorspace})$
says that the Hilbert function of
$I_{(f_1, \ldots, f_s)}$ is
of the form
$(\ref{hilbertfunction})$
and, in addition, that
$I_{(f_1, \ldots, f_s)}$
is minimally generated by
$(\ref{minimally})$.
Thus $($a$)$ and $($b$)$ follow.
Since the lexsegment ideal
with the Hilbert function
$(\ref{hilbertfunction})$
is the universal lexsegment ideal
$(\ref{unilex})$,
one has
$|G((I_{(f_1, \ldots, f_s)})^{\lex})| = s$.
Thus $I_{(f_1, \ldots, f_s)}$
is Gotzmann by $($b$)$.
\end{proof}

If $u$ is a monomial of $A$, then we write
$m(u)$ for the largest integer $j$ for which
$x_j$ divides $u$.
A monomial ideal $I$ of $A$ is called {\em stable}
if, for each monomial $u$ belonging to $G(I)$
and for each $1 \leq i < m(u)$,
one has $(x_i u) / x_{m(u)} \in I$.

\begin{Lemma}
\label{stable}
Let $I$ be a stable ideal of $A$ with
$G(I) = \{ u_1, \ldots, u_s \}$.
Then the Hilbert function of $I$ is
\begin{eqnarray}
\label{stableformula}
H(I,t) = \sum_{j=1}^{s}
{t - \deg u_j + n -m(u_j) \choose n - m(u_j)}.
\end{eqnarray}
\end{Lemma}

\begin{proof}
By virtue of \cite[Lemma 1.1]{EliKer} every monomial
$v \in I$ can be uniquely expressed of the form
$v = uw$, where $u$ and $w$ are monomial of $A$
with $u \in G(I)$ and
$w \in K[x_{m(u)}, \ldots, x_n]$.
It then follows that
\[
I = \bigoplus_{j=1}^{s} u_j K[x_{m(u_j)}, \ldots, x_n].
\]
Thus the desired formula $(\ref{stableformula})$ follows.
\end{proof}

\begin{Lemma}
\label{criticalstable}
A monomial ideal $I$ of $A$ which is both critical
and stable is universal lexsegment.
\end{Lemma}

\begin{proof}
Suppose that $I$ is a critical ideal
of type $(a_1, \ldots, a_s)$.
Let $a_0 = 1$ and $b_i = a_i - a_{i-1} \geq 0$ for
$i = 1, 2, \ldots, s$.
What we must prove is that
$I$ coincides with the lexsegment ideal $(\ref{unilex})$.
Since both $I$ and the lexsegment ideal $(\ref{unilex})$
have the Hilbert function
$(\ref{hilbertfunction})$, our claim follows from
\[
\{ x_1^{b_1+1}, x_1^{b_1}x_2^{b_2+1}, \ldots,
x_1^{b_1} \cdots x_{s-1}^{b_{s-1}} x_{s}^{b_s+1} \}
\subset G(I).
\]

Since
$H(I, a_1 - 1) = 0$ and $H(I, a_1) > 0$
and since $I$ is stable, it follows that
$x_1^{a_1} = x_1^{b_1 + 1} \in G(I)$.
Let $1 \leq p < s$
and
\[
Q =
\{ x_1^{b_1+1}, x_1^{b_1}x_2^{b_2+1}, \ldots,
x_1^{b_1} \cdots x_{p-1}^{b_{p-1}} x_{p}^{b_p+1} \}.
\]
Suppose that $Q \subset G(J)$.
Our work is to show that
$x_1^{b_1} \cdots x_{p}^{b_p}x_{p+1}^{b_{p+1}+1}$
belongs to $G(I)$.
\medskip

We claim that
each monomial $u \in G(I) \setminus Q$
satisfies $m(u) > p$.
Suppose that there is $w \in G(I) \setminus Q$
with $m(w) \leq p$.  Since $I$ is stable,
Lemma \ref{stable} says that
\[
H(I,t)
\geq \sum_{j=1}^{p}
{t - a_j + n - j \choose n - j}
+ {t - \deg w + n - p \choose n - p}.
\]
Since the Hilbert function $H(I, t)$
is equal to $(\ref{hilbertfunction})$,
it follows that
\begin{eqnarray}
\label{binomial}
\sum_{j=p+1}^{s}
{t - a_j + n - j \choose n - j}
\geq {t - \deg w + n - p \choose n - p}
\end{eqnarray}
for all $t$.  However, for $t \gg 0$,
the right-hand side of $(\ref{binomial})$ is
a polynomial on $t$ of degree $n - p$
and the left-hand side of
$(\ref{binomial})$ is that of degree
at most $n - p - 1$.
Hence the inequalities $(\ref{binomial})$
cannot be valid for $t \gg 0$.
This completes the proof of our claim that
each monomial $u \in G(I) \setminus Q$
satisfies $m(u) > p$.

Let $J$ be the universal lexsegment ideal
of $A$ with $G(J) = Q$.  Since
$H(I, t) = H(J, t)$ for $t < a_{p+1}$
and $H(I, a_{p+1}) > H(J, a_{p+1})$,
it follows that there is
a monomial belonging to
$G(I) \setminus Q$
of degree $a_{p+1}$
and that
each monomial of $A$
of degree $a_{p+1}$ belonging to
$I \setminus J$ must belong to $G(I)$.

Let $u = x_1^{c_1} x_2^{c_2} \cdots x_n^{c_n}$
be a monomial belonging to $G(I) \setminus Q$
of degree $a_{p+1}$.
Then
\begin{eqnarray}
\label{power}
c_1 \geq b_1, \ldots, c_{p-1} \geq b_{p-1}.
\end{eqnarray}
To see why $(\ref{power})$ is true,
suppose that there is $1 \leq i < p$ with
$c_i < b_i$.  Since $I$ is stable, one has
\[
v = x_1^{c_1} x_2^{c_2} \cdots x_i^{c_i}
x_{i+1}^{a_{p+1} - \sum_{j=1}^{i}c_j} \in I.
\]
Since $u \not\in J$, it follows that
the monomial
$x_1^{c_1} x_2^{c_2} \cdots x_{i}^{c_{i}}$
can be divided by none of the monomials
$x_1^{b_1} \cdots x_{j-1}^{b_{j-1}} x_{j}^{b_{j}+1}$
with $1 \leq j \leq i$.
Since
$c_i < b_i$, the monomial $v$ can be divided by
none of the monomials
$x_1^{b_1} \cdots x_{j-1}^{b_{j-1}} x_{j}^{b_{j}+1}$
with $i < j \leq p$.
Hence $v \not\in J$.
Since $\deg v = a_{p+1}$, one has $v \in G(I)$.
Since $v \in G(I) \setminus Q$,
the claim stated above
says that $m(v) > p$, a contradiction.
This completes the proof of $(\ref{power})$.

On the other hand, if there is $1 \leq i \leq p$
with $c_i > b_i$, then
$u$ is divided by
$x_1^{b_1} \cdots x_{i-1}^{b_{i-1}} x_{i}^{b_{i}+1}$,
a contradiction.  Hence
\[
c_1 = b_1, \ldots, c_{p-1} = b_{p-1}, c_p \leq b_p.
\]
Since $I$ is stable with $\deg w = a_{p+1}$,
the monomial
\[
w = x_1^{b_1} \cdots x_{p-1}^{b_{p-1}}
x_p^{b_p} x_{p+1}^{b_{p+1}+1}
\]
must belong to $I$.
Since $w \not\in J$, it follows that
$w \in G(I)$, as desired.
\end{proof}

Let $I$ be a homogeneous ideal of $A$.
When $K$ is infinite,
given a monomial order $\sigma$ on $A$,
we write $\gin_\sigma(I)$
for the generic initial ideal
(\cite{Eisenbud} and \cite{Green})
of $I$ with respect to $\sigma$.

\begin{Lemma}
\label{gin}
Suppose that $K$ is infinite.
Let $I$ be a critical ideal of $A$.
Then, for an arbitrary monomial order $\sigma$
on $A$ induced by the ordering
$x_1 > \cdots > x_n$, the generic initial ideal
$\gin_\sigma(I)$ is stable.
Thus in particular
$\gin_\sigma(I)$ is universal lexsegment.
\end{Lemma}

\begin{proof}
Since both $I$ and $\gin_\sigma(I)$
have the same Hilbert function,
it follows that
$\gin_\sigma(I)$ is critical.
Since $\gin_\sigma(I)$ is a monomial ideal,
it follows from \cite[Corollary 1.7]{MuHi} that
$\gin_\sigma(I)$ is Gotzmann.
Thus,
by Gotzmann's persistence theorem (see Lemma \ref{persistence}),
$\gin_\sigma(I)$
is componentwise linear \cite{HerHi}.
Then \cite[Lemma 1.4]{CHH} says that
$\gin_{<_{\rev}}(\gin_\sigma(I))$
$= \gin_\sigma(I)$
is stable.
Here $<_{\rev}$ is the reverse lexicographic
order on $A$ induced by the ordering
$x_1 > \cdots > x_n$.
Since $\gin_\sigma(I)$
is both critical and stable,
it follows from Lemma \ref{criticalstable}
that $\gin_\sigma(I)$ is universal lexsegment.
\end{proof}

Note that if $\mathrm{char}(K)=0$ then Lemma \ref{gin} is obvious
since generic initial ideals are stable in characteristic $0$.

\begin{Lemma}
\label{crucial}
Suppose that a homogeneous ideal $I$ of $A$
is a critical ideal of type $(a_1, \ldots, a_s)$,
where $2 \leq s \leq n$.  Then there exists
a homogeneous polynomial $f$ of $A$ with
$\deg f = a_1 - 1$ together with
a homogeneous ideal $J$ of $A$ such that
\[
I = f \cdot J.
\]
\end{Lemma}

\begin{proof}
To prove the statement, by considering an extension field,
we may assume
that the base field $K$ is infinite.
Then there is a linear transformation
$\varphi$ with
$\initial_{<_{\lex}}(\varphi(I)) = \gin_{<_{\lex}}(I)$.
Considering $\varphi(I)$ instead of $I$,
one may assume that
$\initial_{<_{\lex}}(I) = \gin_{<_{\lex}}(I)$.
Lemma \ref{gin} says that
$\initial_{<_{\lex}}(I)$ is universal lexsegment.
Hence
\[
\initial_{<_{\lex}}(I)
= (x_1^{b_1+1}, x_1^{b_1}x_2^{b_2+1}, \ldots,
x_1^{b_1} \cdots x_{s-1}^{b_{s-1}} x_s^{b_s+1}),
\]
where $b_i = a_i - a_{i-1}$,
$1 \leq i \leq s$,
with $a_0 = 1$.

To simplify the notation,
let $u_i = x_1^{b_1} \cdots x_i^{b_i}$
for $i = 1, \ldots, s$.  Thus
$\initial_{<_{\lex}}(I)
= (u_1 x_1, \ldots, u_s x_s)$.
Let ${\mathcal G} = \{ g_1, \ldots, g_s \}$
be a Gr\"obner basis of $I$,
where $g_i$ is a homogeneous polynomial
of $A$ with
$\initial_{<_{\lex}}(g_i) = u_i x_i$
for each $1 \leq i \leq s$,
and ${\mathcal G}' = \{ g_2, \ldots, g_s \}$.
We show that ${\mathcal G}'$ is
a Gr\"obner basis with respect to
$<_{\lex}$.
Let
$2 \leq i < j \leq s$
and divide the $S$-polynomial of $g_i$ and $g_j$
by ${\mathcal G}$, say,
\[
(u_j / u_i) x_j g_i - x_i g_j
= p_1 g_1 + \cdots + p_s g_s,
\]
where $p_1, \ldots, p_s$ are homogeneous
polynomials of $A$ with
\[
\initial_{<_{\lex}}(p_k g_k)
\leq_{\lex}
\initial_{<_{\lex}}((u_j / u_i) x_j g_i - x_i g_j)
\]
for each $1 \leq k \leq s$.
Since
\[
\initial_{<_{\lex}}(p_1g_1)
= \initial_{<_{\lex}}(p_1) x_1^{b_1 + 1}
<_{\lex}
\initial_{<_{\lex}}(x_ig_j)
= x_i x_j u_j
= x_i x_j (x_1^{b_1} \cdots x_{j}^{b_j}),
\]
it follows that $p_1 = 0$.
In other words,
a remainder of
the $S$-polynomial of $g_i$ and $g_j$
with respect to ${\mathcal G}'$
can be $0$.
Hence ${\mathcal G}'$ is a Gr\"obner basis
with respect to $<_{\lex}$, as desired.

Now, we prove Lemma \ref{crucial}
by using induction on $s$.
Let $s = 2$.
Again, divide the $S$-polynomial of $g_1$ and $g_2$
by ${\mathcal G}$, say,
\[
x_2^{b_2 + 1} g_1 - x_1 g_2
= p_1 g_1 + p_2 g_2,
\]
where $p_1$ and $p_2$
are homogeneous with
$\deg p_1 = b_2 + 1$
and $\deg p_2 = 1$,
and where each of $p_1$ and $p_2$
satisfies
\[
\initial_{<_{\lex}}(p_k g_k)
\leq_{\lex}
\initial_{<_{\lex}}(x_2^{b_2 + 1} g_1 - x_1 g_2).
\]
One has
$(x_2^{b_2 + 1} - p_1)g_1
= (x_1 + p_2) g_2$.
Since
$\initial_{<_{\lex}}(x_2^{b_2 + 1} - p_1)
= x_2^{b_2 + 1}$ and since
$\initial_{<_{\lex}}(x_1 + p_2)
= x_1$, the polynomial $x_1 + p_2$
must divide $g_1$.
Let $f = g_1 / (x_1 + p_2)$.
Then $\deg f = a_1 - 1$ and
$I =
f \cdot (x_1 + p_2, x_2^{b_2 + 1} - p_1)$,
as required.

Next, let $s > 2$ and write $J$ for the ideal of $A$
generated by ${\mathcal G}'$.  Since
\[
\initial_{<_{\lex}}(J)
= (u_2 x_2, \ldots, u_s x_s)
= x_1^{b_1}
(x_2^{b_2 + 1}, x_2^{b_2} x_3^{b_3 + 1},
\ldots, x_2^{b_2} \cdots x_{s-1}^{b_{s-1}}
x_s^{b_s + 1})
\]
and since
\[
(x_2^{b_2 + 1}, x_2^{b_2} x_3^{b_3 + 1},
\ldots, x_2^{b_2} \cdots x_{s-1}^{b_{s-1}}
x_s^{b_s + 1})
\]
is universal lexsegment
by a permutation of the variables,
the ideal $J$ is
a critical ideal of type
$(a_2, \ldots, a_s)$.
The induction hypothesis guarantees
the existence of a homogeneous polynomial
$f_0$ of $A$ with $\deg(f_0) = a_2 - 1$
which divides each of
$g_2, \ldots, g_s$.
Since $\initial_{<_{\lex}}(f_0)$ divides
$\initial_{<_{\lex}}(g_i) = u_i x_i$
for each $1 < i \leq s$,
one has
$\initial_{<_{\lex}}(f_0) = u_2$.
Let $g_i' = g_i / f_0$
for $i = 2, \ldots, s$.
Thus in particular
$\initial_{<_{\lex}}(g_2')
=u_2x_2 /u_2= x_2$.

Now, divide the $S$-polynomial of $g_1$ and $g_2$
by ${\mathcal G}$, say,
\[
x_2^{b_2 + 1} g_1 - x_1 (f_0 g_2')
= q_1 g_1 + q_2 (f_0 g_2')
+ \cdots + q_s (f_0 g_s'),
\]
where $q_1, \ldots, q_s$ are homogeneous
polynomials of $A$ with
\[
\initial_{<_{\lex}}(q_1 g_1)
\leq_{\lex}
\initial_{<_{\lex}}(x_2^{b_2 + 1} g_1 - x_1 (f_0 g_2'))
\]
and with
\[
\initial_{<_{\lex}}(q_k (f_0 g_k'))
\leq_{\lex}
\initial_{<_{\lex}}(x_2^{b_2 + 1} g_1 - x_1 (f_0 g_2'))
\]
for each $2 \leq k \leq s$.
Let
\[
f_0 h
= q_2 (f_0 g_2')
+ \cdots + q_s (f_0 g_s').
\]
Thus
\[
(x_2^{b_2 + 1} - q_1)g_1
= f_0 (x_1 g_2' + h).
\]
Since
$\initial_{<_{\lex}}
(x_2^{b_2 + 1} - q_1)
= x_2^{b_2 + 1}$,
$\initial_{<_{\lex}}
(g_1) = x_1^{b_1 + 1}$
and
$\initial_{<_{\lex}}
(x_1 g_2' + h)
= x_1 x_2$,
it follows that
$x_1 g_2' + h$
can divide neither
$x_2^{b_2 + 1} - q_1$
nor $g_1$.
Thus
$x_1 g_2' + h$ is a product
$(x_1 + h_1)(x_2 + h_2)$,
where $h_1$ and $h_2$ are
homogeneous polynomials of $A$
with $\deg h_1 = \deg h_2 = 1$,
such that
$x_1 + h_1$ divides $g_1$
and $x_2 + h_2$ divides
$x_2^{b_2 + 1} - q_1$.
Let $f = g_1 / (x_1 + h_1)$.
Then $\deg f = a_1 - 1$
and $f$ divides
both $g_1$ and $f_0$.
\end{proof}

\section{Proof of Theorem \ref{gotzmann}}
We are now in the position to give a proof of
Theorem \ref{gotzmann}.

\begin{proof}[Proof of Theorem \ref{gotzmann}]
First, Corollary \ref{Sydney}
guarantees $($ii$) \Rightarrow ($iii$)$.
Second,
$($iii$) \Rightarrow ($i$)$ is clear by the definition of Gotzmann ideals
and that of critical functions.

On the other hand, a proof of
$($i$) \Rightarrow ($ii$)$ will be achieved
by induction on $s$.
Let $I \subset A$ be a critical ideal
of type $(a_1, \ldots, a_s)$.
Let $s = 1$.  Let $f_1$ be a homogeneous
polynomial of degree $a_1$ belonging to $I$.
Then the Hilbert function of the ideal
$( f_1 )$ of $A$ coincides with
that of $I$.  Thus $I = ( f_1 )$,
as desired.

Let $s > 1$.
Lemma \ref{crucial} guarantees that
$I = f \cdot J$, where
$f$ is a homogeneous polynomial of $A$
with $\deg f = a_1 - 1$ and where $J$ is
a homogeneous ideal of $A$.
The Hilbert function of $J$ is
$H(J, t) = H(I, t + a_1 - 1)$.
Hence $J$ is a critical ideal of
type$(1, a_2 - a_1 + 1, \ldots, a_s - a_1 + 1)$.
Since $H(J,1) \neq 0$,
there exists a linear transformation
$\varphi$ on $A$ with
$x_1 \in \varphi(J)$.  Let $J'$ be
the ideal
\[
J' = \varphi(J) \cap K[x_2, \ldots, x_n]
\]
of $K[x_2, \ldots, x_n]$.  Then
\[
\varphi(J) = x_1 K[x_1, \ldots, x_n] \bigoplus J'.
\]
Since the Hilbert function of $J'$ is
\[
H(J',t) = H(\varphi(J), t) - H(x_1 K[x_1, \ldots, x_n], t),
\]
the ideal $J'$ of
$K[x_2, \ldots, x_n]$ is a critical ideal of
type $(a_2 - a_1 + 1, \ldots, a_s - a_1 + 1)$.
The induction hypothesis then guarantees
the existence of a linear transformation
$\psi$ on $K[x_2, \ldots, x_n]$
such that
$\psi(J')$ is a canonical critical ideal
of $K[x_2, \ldots, x_n]$,
say
\[
\psi(J') = (f_2x_2, \ldots,
f_2 \cdots f_{s-1} x_{s-1},
f_2 \cdots f_s),
\]
where
$f_i \in K[x_i, x_{i+1}, \ldots, x_n]$
for each $2 \leq i \leq s$
and where $\deg f_s > 0$.
Now,
regarding $\psi$
to be a linear transformation
on $A$ by setting $\psi(x_1) = x_1$,
one has
\begin{eqnarray*}
(\psi \circ \varphi)(I) & = &
((\psi \circ \varphi)(f)) \cdot ((\psi \circ \varphi)(J)) \\
& = &
((\psi \circ \varphi)(f)) \cdot (\psi(x_1 A \bigoplus J')) \\
& = &
(\psi \circ \varphi)(f) \cdot (x_1 A \bigoplus \psi(J')).
\end{eqnarray*}
Let $f_1 = (\psi \circ \varphi)(f)$.
Then it follows that
\[
(\psi \circ \varphi)(I)
= (f_1x_1, f_1f_2x_2, \ldots,
f_1 f_2 \cdots f_{s-1} x_{s-1},
f_1 f_2 \cdots f_s)
\]
as desired.
\end{proof}

\begin{Remark}
{\em
In Lemma \ref{gin} we assume that the base field $K$ is infinite.
However, this assumption is not required in Theorem \ref{gotzmann}
since Lemma \ref{gin} is only required to prove Lemma \ref{crucial}
and since we may assume that $K$ is infinite to prove Lemma \ref{crucial}.
}
\end{Remark}

\begin{Example}
{\em
For a permutation $\pi$ on
$[n] = \{ 1, \ldots, n \}$
and for a monomial
$u = x_1^{a_1} \cdots x_n^{a_n}$
of $A$, we set
$\pi(u) =
x_{\pi(1)}^{a_1} \cdots x_{\pi(n)}^{a_n}$.
A monomial ideal $I$ of $A$
with $G(I) = \{ u_1, \ldots, u_s \}$ is
called {\em trivially Gotzmann}
if there is a permutation $\pi$ on
$[n]$ such that the monomial ideal
$(\pi(u_1), \ldots, \pi(u_s))$
is lexsegment.

Given an arbitrary integer $m > n$,
there exists a Gotzmann monomial ideal $I$
of $A$ with $|G(I)| = m$
which cannot be trivially Gotzmann.
In fact, let
$q \geq 2$
be an integer and
$I$ the monomial ideal of $A$
generated by
$x_1^{q-1}x_2, x_1^{q-2}x_2^2,
\ldots, x_1 x_2^{q-1}$
together with
$x_3, \ldots, x_n$.
Then $I$ is a Gotzmann ideal
with $|G(I)| = n + q - 3$
which cannot be trivially Gotzmann.
}
\end{Example}

\section{Inflexible Hilbert functions}

Throughout this section,
we assume that $A=K[x_1,\dots,x_n]$ is the polynomial ring over
an infinite field $K$ with each $\deg x_i =1$.
The \textit{graded Betti numbers} of a homogeneous ideal $I$
are the integers $\beta_{i,j}(I) = \dim_K \Tor_i (I,K)_j$.
In other words, the graded Betti numbers 
appear in the 
minimal graded free resolution 
\[
\ \ \ \ \ \cdots
\longrightarrow
\bigoplus_{j} S(-j)^{\beta_{2,j}(I)}
\longrightarrow
\bigoplus_{j} S(-j)^{\beta_{1,j}(I)}
\longrightarrow
\bigoplus_{j} S(-j)^{\beta_{0,j}(I)}
\longrightarrow I
\longrightarrow 0
\]
of $I$ over $A$. 
By the famous result of Bigatti \cite{Big},
Hulett \cite{Hul} and Pardue \cite{Par},
there exists an ideal having the largest graded Betti numbers
among all ideals for a fixed Hilbert function.
However, an ideal having the smallest graded Betti numbers among
all ideals for a fixed Hilbert function
need not exist.
Recently, existence and non-existence of the smallest graded Betti numbers of
ideals for a fixed Hilbert function are studied in several papers
(see e.g., \cite{CE,DMMR,R}).
An extremal example of a Hilbert function having an ideal with
the smallest graded Betti numbers is a Hilbert function $H$
for which all ideals with the Hilbert function $H$ have
the same graded Betti numbers.
Let $H : \ZZ_{\geq 0} \to \ZZ_{\geq 0}$ be
the Hilbert function of a homogeneous ideal of $A$.
We say that $H$ is an \textit{inflexible Hilbert function of $A$}
if all homogeneous ideals of $A$ with the Hilbert function $H$
have the same graded Betti numbers.
It is known \cite[Corollary 1.4]{HerHi} that
$H$ is an inflexible Hilbert function of $A$
if and only if
all ideals of $A$ with the Hilbert function $H$ are Gotzmann.
Thus, in particular,
Theorem \ref{gotzmann} shows that critical functions are 
inflexible Hilbert functions.
In this section,
we introduce some more inflexible Hilbert functions.

First, we recall fundamental properties on Macaulay
representations and the minimal growth of Hilbert functions.
Given positive integers $a$ and $d$,
there exists the unique representation of $a$,
called the {\em $d$-th Macaulay representation}
of $a$, of the form
\begin{eqnarray}
\label{Macaulay}
a = {a(d) + d  \choose d}
+ \cdots +
{a(k) + k \choose k},
\end{eqnarray}
where $ k \geq 1$ and where
$a(d) \geq \cdots \geq a(k) \geq 0$.
We recall the following easy fact
(see \cite[Lemma 4.2.7]{BruHer}).

\begin{Lemma}
\label{macaulayexpand}
Let $a_d \geq \cdots \geq a_1 \geq -1$ and
$b_d  \geq \cdots \geq b_1 \geq -1$ be integers.
Then
$${a_d + d  \choose d}
+ \cdots + {a_1 + 1 \choose 1} >
{b_d + d  \choose d}
+ \cdots +
{b_1 + 1 \choose 1}$$
if and only if
$(a_d,\dots,a_1) >_{\mathrm{lex}} (b_d,\dots,b_1)$.
\end{Lemma}

By using the $d$-th Macaulay representation
$(\ref{Macaulay})$
of a positive integer $a$, one defines
\[
a \Macd
= {a(d) + d + 1 \choose d}
+ \cdots +
{a(k) + k + 1 \choose k}
\]
and
$0 \Macd =0$.
This number is convenient to describe the minimal growth of Hilbert functions
of homogeneous ideals.
Let $L \subset A_d$ be a $K$-vector space spanned by a lexsegment set
of monomials of degree $d$.
Let $u$ be the minimal monomial in $L$ with respect to $<_{\lex}$.
Then $u$ can be written in the form
$$u= x_1^{a_0-1} x_2 ^{a_1-a_0} \cdots x_k^{a_{k-1} -a_{k-2}} x_{k+1}^{a_k -a_{k-1}} x_n^{d-a_k},$$
where $0 <a_0 \leq a_1 \leq \cdots \leq a_k \leq d$
and $0 \leq k \leq n-2$.
Then it follows from \cite[Lemma C.10]{IK} that
$$\dim_K L=
{d-a_0 +n-1 \choose n-1} + \cdots + { d-a_k + n-1-k \choose n-1-k}$$
is the $(n-1)$-th Macaulay representation of $\dim_K L$.
On the other hand, since $A_1 \cdot L$ is also lexsegment 
and since the minimal element of $A_1 \cdot L$ w.r.t.\ $<_{\lex}$ is
$x_n u$,
one has
\begin{eqnarray}
\label{MacEq}
 \dim_K (A_1 \cdot L)= (\dim_K L) \Mac.
\end{eqnarray}
Hence by (\ref{MacIneqSpace}) any homogeneous ideal $I$ of $A$ satisfies
$$
H(I,k+1) \geq H(I,k) \Mac
\ \ \mbox{ for all }k \geq 0.
$$

Suppose that $H$ is the Hilbert function of a homogeneous ideal of $A$.
Peeva \cite[Corollary 1.4]{Peeva} proved that,
any numerical function $H$
satisfying that,
for all $k \geq 1$,
if $H(k) \Mac <H(k+1)$, then $H(k+1) \Mac = H(k+2)$,
is an inflexible Hilbert function.
If $n=2$, those functions essentially characterize 
inflexible Hilbert functions of $A$.
Indeed, it is not hard to see that,
in case of $n=2$, $H$ is an inflexible Hilbert function of $A$
if and only if $H$ satisfies that,
for all $k \geq 1$,
if $H(k) \Mac <H(k+1)$ and $H(k) \ne 1$, then $H(k+1) \Mac = H(k+2)$.
See \cite{I} for further results on graded Betti numbers
of homogeneous ideals of $K[x_1,x_2]$ for a fixed Hilbert function.
Peeva's result and Theorem \ref{gotzmann} lead us to consider
the following numerical functions.

\begin{Definition}
{\em
\label{segmentwise}
A numerical function $H : \mathbb{Z}_{\geq 0} \to \mathbb{Z}_{\geq 0}$ is said to be \textit{segmentwise critical} if there exist integers
$0=s_0 < s_1< \cdots < s_\ell = \infty$ such that
\begin{itemize}
\item[(i)]
for each $j=0,1,\dots,\ell-1$, 
there exists a critical function $H_j$ with 
$H(t)=H_j(t)$ for all $s_j \leq t \leq s_{j+1}$;
\item[(ii)]
$H(s_{j+1}-1) \Mac =H(s_{j+1})$ for $j=0,1,\dots,\ell-2$.
\end{itemize}
}
\end{Definition}

Let $d \geq 0$ be an integer.
We say that a subspace $V$ of $A_d$ is a \textit{critical space}
if there exists a critical ideal $I$ of $A$ such that
$V=I_d$.
Note that critical spaces are Gotzmann by Theorem \ref{gotzmann}.
For any homogeneous ideal $I$ of $A$,
we write $\gin (I) = \gin_{<_{\rev}}(I)$.
The goal of this section is to prove the following result.

\begin{Theorem}
\label{2-5}
If $I$ has a segmentwise critical Hilbert function then
$I_k$ is a critical space for all $k \geq 0$.
In particular, $I$ is Gotzmann and $\gin(I)$ is lexsegment.
\end{Theorem}

To prove Theorem \ref{2-5},
we first study critical spaces.
For a homogeneous ideal $I$ of $A$, the ideal
$$ I^\sat= 
\bigcup_{k \geq 0} (I: (x_1,\dots,x_n)^k)$$
is called the \textit{saturation of $I$}.
The \textit{(Castelnuovo--Mumford) regularity} of a homogeneous ideal
$I$ is the integer 
$$\reg(I)=\max\{ j: \beta_{i,i+j}(I) \ne 0 \mbox{ for
some }i\}.$$
The following facts are known
(see \cite[\S 2]{Green}).

\begin{Lemma}
\label{bayerstillman}
Let $I$ be a homogeneous ideal of $A$. 
\begin{itemize}
\item[(i)]
$I^\sat_k = I_k$ for all $k \geq \reg(I)$;
\item[(ii)] $ \gin(I^\sat) 
= \bigcup_{k \geq 0} (\gin(I): x_n^k) =(\gin(I))^\sat$;
\item[(iii)] (Bayer--Stillman)
$\reg (I) = \reg(\gin(I))$.
\end{itemize}
\end{Lemma}

Note that the second equation of (ii)
follows from \cite[Proposition 15.24]{Eisenbud}.

A homogeneous ideal $I$ is said to have a \textit{$d$-linear resolution}
if it is generated in degree $d$ and $\reg(I)=d$.
The following result of Gotzmann is called Gotzmann's persistence
theorem (see \cite[\S 4.3]{BruHer} and \cite{Green} for the proof).

\begin{Lemma}[Gotzmann's persistence theorem \cite{Gpersist}]
\label{persistence}
Let $V \subset A_d$ be a Gotzmann space and $I$ the ideal generated by
$V$. Then \begin{itemize}
\item[(i)] $A_1 \cdot V$ is Gotzmann;
\item[(ii)] $I$ has a $d$-linear resolution.
\end{itemize}
\end{Lemma}

Recall that a homogeneous ideal $I \subset A$ is
{\em saturated} if $I : (x_1, \ldots, x_n) = I$.
Note that a homogeneous ideal $I$ of $A$ is saturated if and only if
the depth of $A/I$ is positive (see \cite[\S 1.2]{BruHer}).
The next result immediately follows from \cite[Corollary 1.4]{MuHi}.

\begin{Lemma}
\label{satlex}
A lexsegment ideal $I$ of $A$ is saturated if and only if
$|G(I)| \leq n-1$.
\end{Lemma}

Let $V$ be a subspace of $A_d$ and $I$ the ideal of $A$
generated by $V$.
Let $\gin(V)$ be the \textit{generic initial space of $V$},
in other words,
$$\gin(V) = \gin(I)_d.$$

\begin{Lemma}
\label{criticallex}
Let $V$ be a subspace of $A_d$ and $I$ the ideal of $A$ generated by $V$. 
If $I$ has a $d$-linear resolution and $\gin(V)=\Lex(V)$, then
\begin{itemize}
\item[(i)] 
$I^\sat$ is a critical ideal generated by at most $n-1$
homogeneous polynomials;
\item[(ii)]
$V$ is a critical space.
\end{itemize}
\end{Lemma}

\begin{proof}
By the assumption and Lemma \ref{bayerstillman} (iii),
it follows that $\gin(I)$ has a $d$-linear resolution.
Thus $\gin(I)$ is generated by $\gin(I)_d= \gin(V)=\Lex(V)$.
Hence $\gin(I)$ is lexsegment.
It is clear that the saturation of a lexsegment ideal
is again a lexsegment ideal.
Thus, by Lemma \ref{satlex}, the ideal $\gin(I)^\sat$ is a universal
lexsegment ideal with $|G(\gin(I)^\sat)| \leq n-1$.
Since $\gin(I^\sat)=\gin(I)^\sat$ by Lemma \ref{bayerstillman} (ii),
the Hilbert function of $I^\sat$ is equal to that of the universal lexsegment ideal $\gin(I)^\sat$.
Hence $I^\sat$ is a critical ideal.
Also since the number of generators of $I^\sat$ is smaller than
or equal to that of $\gin(I^\sat)=\gin(I)^\sat$,
the ideal $I^\sat$ is
generated by at most $n-1$ homogeneous polynomials.
Finally,
since $\reg(I)=d$, we have
$V=I_d = (I^\sat)_d$ by Lemma \ref{bayerstillman} (i),
and therefore
$V$ is a critical space.
\end{proof}

Lemma \ref{criticallex} 
enable us to characterize critical spaces
in terms of generic initial spaces.

\begin{Corollary}
\label{classify}
A subspace $V$ of $A_d$ is a critical space if and only if
$V$ is Gotzmann and $\gin(V)=\Lex(V)$.
\end{Corollary}

\begin{proof}
If $V$ is a critical space then
$V$ is Gotzmann and $\gin(V)=\Lex(V)$ by Theorem \ref{gotzmann} and Lemma \ref{gin}.
On the other hand,
by Lemmas \ref{persistence} and \ref{criticallex},
if $V$ is Gotzmann and $\gin(V)=\Lex(V)$, then 
$V$ is a critical space. 
\end{proof}

\begin{Example}
{\em
Let $A=K[x_1,x_2,x_3]$,
$I=(x_1^3,x_1^2x_2,x_1x_2^2,x_2^3,x_1^2x_3,x_1x_2x_3,x_2^2x_3)$
and $V=I_3$.
It is easy to see that $V$ is a Gotzmann space.
However, $V$ is not a critical space.
Indeed, $\gin(V)=V \ne \Lex(V)$ and $I^\sat=(x_1^2,x_1x_2,x_2^2)$ does not have a critical Hilbert function.
On the other hand,
for any $f \in A$,
the ideal 
$J=(fx_1^2,f x_1x_2,fx_1x_3,f(x_2+x_3)x_2,f(x_2+x_3)x_3)$
is generated by a critical space. Indeed,
 $J^\sat =(fx_1,f(x_2+x_3))$
is a canonical critical ideal and $J$ is generated by the homogeneous component of $J^\sat$ of degree $\deg f +2$.
Note also that these critical spaces are parametrized by some Hilbert schemes
(see the appendix).
}
\end{Example}

By Lemma \ref{criticallex},
if $V \subset A_d$ is a critical space,
then the saturation of the ideal $I$ generated by $V$
is critical.
We know the type of this critical ideal
by using the Macaulay representation of $\dim_K V$.

\begin{Lemma}
\label{2-3}
Let $V \subset A_d$ be a critical space
and $I$ the ideal of $A$ generated by $V$.
Let
$\dim_K V = \sum_{j=1}^p {d-a_p + n-j \choose n-j}$
be the $(n-1)$-th Macaulay representation of $\dim_K V$.
Then $I^\sat$ is a critical ideal of type $(a_1,\dots,a_p)$.
\end{Lemma}

\begin{proof}
Suppose that $I^\sat$ is a critical ideal of type $(b_1,\dots,b_q)$.
What we must prove is $(a_1,\dots,a_p)=(b_1,\dots,b_q)$.

By Theorem \ref{gotzmann} and Corollary \ref{Sydney},
the type of a critical ideal is equal to the sequence
of degrees of its minimal generators.
Thus, by Lemma \ref{criticallex}, we have $q \leq n-1$.
In particular,
$$ 
H(I^\sat,k)= \sum_{j=1}^q { k-b_j + n-j \choose n-j}
$$
is the $(n-1)$-th Macaulay representation of $H(I^\sat,k)$ for $k \geq b_q$.
On the other hand, (\ref{MacEq}) and Gotzmann's persistence theorem imply
$$
H(I,k)= \sum_{j=1}^p { k-a_j + n-j \choose n-j}
\ \ \ \ \mbox{ for all }k \geq d.
$$
Since the above equation is the $(n-1)$-th Macaulay representation
of $H(I,k)$ for all $k \geq d$
and since $H(I,k)=H(I^\sat,k)$ for $k \gg 0$,
it follows from Lemma \ref{macaulayexpand} that $(a_1,\dots,a_p)=(b_1,\dots,b_q)$.
\end{proof}

Next,
we study the relation between critical functions
and Macaulay representations.
Let $H$ be the critical function of type $(a_1,\dots,a_p)$.
If $p \leq n-1$ then (\ref{hilbertfunction}) coincides with the $(n-1)$-th
Macaulay representation.
On the other hand, in case of $p=n$,
the $(n-1)$-th Macaulay representation of $H(k)$ is given by the following formula.

\begin{Lemma}
\label{2-0}
Let $a_0=0<a_1  \leq \cdots \leq a_n \leq d$ be integers
and let $a= \sum_{j=1}^n { d-a_j + n-j \choose n-j}$.
Set $s= \min \{ k : a_k = a_{n-1},\ k \leq n-1\}$.
Then
$$a= \left\{ \sum_{j=1}^{s-1} { d-a_j + n-j \choose n-j} \right\}
+{t-(a_s-1)+n-s \choose n-s}$$
is the $(n-1)$-th Macaulay representation of $a$.
\end{Lemma}

\begin{proof}
The statement immediately follows from the next equation
$${t-a_s +n-s \choose n-s } +\cdots +{t-a_s + 1 \choose 1}
+ {t-a_n +0 \choose 0} = { t-(a_s-1) + n-s \choose n-s}.$$
\end{proof}

For  ${\bf a} = (a_1,\dots,a_p) \in \mathbb{Z}_{>0}^p$, where $0<a_1 \leq \cdots \leq a_p$
and where $1 \leq p \leq n$, write
\begin{eqnarray*}
\bar {\bf a} = 
\left\{ \begin{array}{ll}
(a_1,\dots,a_p), & \mbox{ if } p<n,\\
(a_1,\dots,a_{s-1},a_s-1), & \mbox{ if } p=n,\\
\end{array}
\right.
\end{eqnarray*}
where $s= \min \{ k : a_k = a_{n-1},\ k \leq n-1\}$.
We also require the following fact, which immediately follows
from Gotzmann's persistence theorem.

\begin{Lemma}
\label{obvious}
Let $I$ be a homogeneous ideal of $A$.
If $H(I,k) \Mac = H(I,k+1)$, then
$I_k$ and $I_{k+1}$ are Gotzmann spaces.
\end{Lemma}

\begin{proof}
Clearly $I_k$ is Gotzmann
and $I$ has no generators of degree $k+1$.
Hence $I_{k+1} = A_1 \cdot I_k$ is also a Gotzmann space
by Gotzmann's persistence theorem.
\end{proof}

Now we prove Theorem \ref{2-5}.

\begin{proof}[Proof of Theorem \ref{2-5}]
For any integer $k \geq 0$,
we write $I_{\langle k \rangle}$ for
the ideal generated by all polynomials in $I$ of degree $k$.
Let $0=s_0 < s_1 < \cdots < s_\ell = \infty$ be integers
satisfying the conditions of Definition \ref{segmentwise}.
Let
$$I[i] = (I_{\langle s_i \rangle})^\sat + I_{\langle s_{i} +1 \rangle} +\cdots + I_{\langle s_{i+1} \rangle}
\ \ \ \mbox{ for } i=0,1,\dots,\ell-1.$$
We show that $I[i]$ is a critical ideal
with $(I[i])_k = I_k$ for all $s_i \leq k \leq s_{i+1}$.
Since $I_0$ is a critical space,
by using induction on $i$,
it is enough to show that if $I_{s_i}$ is a critical space,
then $I[i]$ is a critical ideal with 
$(I[i])_k = I_k$ for all $s_i \leq k \leq s_{i+1}$.

First, we show
$(I[i])_k = I_k$ for all $s_i \leq k \leq s_{i+1}$.
Since $I_{s_i}$ is a critical space by the assumption,
$I_{s_i}$ is a Gotzmann space.
Thus, by Lemma \ref{persistence} (ii),
we have $\reg(I_{\langle s_i \rangle})=s_i$.
Hence, by Lemma \ref{bayerstillman} (i), we have
\begin{eqnarray}
\label{hoshi1}
((I_{\langle s_i \rangle})^\sat)_k 
= (I_{\langle s_i \rangle})_k \subset I_k
\ \ \mbox{ for all }k \geq s_i.
\end{eqnarray}
Then it is clear that $(I[i])_k = I_k$ for all $s_i \leq k \leq s_{i+1}$.

Next, we show that $I[i]$ has a critical Hilbert function.
Suppose that, for $k=s_i,s_i+1,\dots,s_{i+1}$,
$H(I,k)$ is equal to the critical function of type 
${\bf a}=(a_1,\dots,a_p)$.
Thus
\begin{eqnarray}
\label{hoshi2}
H(I[i],k)=H(I,k) = \sum_{j=1}^p { k-a_j + n-j \choose n-j}
\ \ \mbox { for } s_i \leq k \leq s_{i+1}.
\end{eqnarray}
We may assume $a_p \leq s_{i+1}$
since ${ k -a_p + n-p \choose n-p}=0$ if $k < a_p$.
Set $\bar {\bf a}=(b_1,\dots,b_q)$.

We claim that
\begin{eqnarray}
\label{hoshi3}
&&H(I[i],k)
=\sum_{j=1}^q { k-b_j + n-j \choose n-j}
=\sum_{j=1}^p { k-a_j + n-j \choose n-j}
\ \ \mbox{ for all $k \geq s_{i+1}$.}
\end{eqnarray}
If $i=\ell-1$ then there is nothing to prove.
Suppose $i<\ell-1$.
Since $a_p \leq s_{i+1}$,
it follows from Lemma \ref{2-0} and (\ref{hoshi2}) that
$$\dim_K I_{s_{i+1}}=H(I,s_{i+1})
= \sum_{j=1}^q {s_{i+1} -b_j + n-j \choose n-j}$$
is the $(n-1)$-th Macaulay representation of
$\dim_K I_{s_{i+1}}$.
On the other hand, by Lemma \ref{obvious} and condition (ii) of Definition \ref{segmentwise},
it follows that $I_{s_{i+1}}$ is a Gotzmann space.
Hence by Gotzmann's persistence theorem we have
\begin{eqnarray*}
&&H(I_{\langle s_{i+1} \rangle},k) 
=\sum_{j=1}^q { k-b_j + n-j \choose n-j}
\ \ \mbox{ for all $k \geq s_{i+1}$.}
\end{eqnarray*}
Then the desired equation (\ref{hoshi3}) follows from Lemma \ref{2-0},
since (\ref{hoshi1}) says that $I[i]_k = (I_{\langle s_{i+1} \rangle})_k$
for $k \geq s_{i+1}$.

Now we show that $I[i]$ has a critical Hilbert function.
\medskip

\textit{Case 1}:
Suppose $a_p \leq s_i$.
We show that $I[i]$ has the critical Hilbert function of type $(b_1,\dots,b_q)$.
Since $a_p \leq s_i$,
Lemma \ref{2-0} together with (\ref{hoshi2}) and (\ref{hoshi3}) implies
\begin{eqnarray}
\label{hoshi4}
H(I[i],k)= \sum_{j=1}^q {k-b_j + n-j \choose n-j}
\ \ \ \ \mbox{ for all }k \geq s_i.
\end{eqnarray}
In particular,
$$\dim_K I_{s_i} = H(I[i],s_i) =\sum_{j=1}^q {s_i-b_j + n-j \choose n-j}$$
is the $(n-1)$-th Macaulay
representation of $\dim_K I_{s_i}$.
Thus, by Lemma \ref{2-3},
one has
$$
H(I[i],k)=H((I_{\langle s_i \rangle})^\sat,k)
= \sum_{j=1}^q {k-b_j + n-j \choose n-j}
\ \ \ \ \mbox{ for }k \leq s_i.
$$
Then the above equation together with (\ref{hoshi4}) implies that
the Hilbert function of $I[i]$
is the critical function of type $(b_1,\dots,b_q)$.
\medskip

\textit{Case 2}:
Suppose $a_p > s_i$.
We show that $I[i]$ has the critical Hilbert function of type $(a_1,\dots,a_p)$.
Let $\ell = \min \{ k: a_k > s_i\} -1$.
Then, by (\ref{hoshi2}),
$$\dim_K I_{s_i} = H(I,s_i)
= \sum_{j=1}^\ell {s_i -a_j + n-j \choose n-j}
$$
is the $(n-1)$-th Macaulay representation of
$\dim_K I_{s_i}$.
Thus, by Lemma \ref{2-3}, we have
$$
H(I[i],k)=H((I_{\langle s_i \rangle})^\sat,k)
= \sum_{j=1}^\ell {k-a_j + n-j \choose n-j}
\ \ \ \ \mbox{ for }k \leq s_i.
$$
Then the above equation together with (\ref{hoshi2}) and
(\ref{hoshi3}) says that
the Hilbert function of $I[i]$ is the critical function
of type $(a_1,\dots,a_p)$.
\end{proof}

\begin{Example}
{\em
Unfortunately, we cannot expect a simple canonical formula like Theorem
\ref{gotzmann} for homogeneous ideals having a segmentwise critical
Hilbert function.
If $I$ has a segmentwise critical Hilbert function,
then each $I_k$ is a critical space by Theorem \ref{2-5}.
However, there does not necessarily exist a linear transformation $\varphi$ on $S$
such that $\varphi(I)_k$
is a canonical critical space for all $k$.
Let
$$I=(
x_1^3,x_1^2x_2,x_1^2x_3^2,x_1(x_2^4+x_3^4)x_3) \subset \QQ[x_1,x_2,x_3].
$$
Let $J=(x_1^3,x_1^2x_2,x_1^2x_3^2)$ and $J'=(x_1^2,x_1(x_2^4+x_3^4)x_3).$
Then $J$ and $J'$ are critical ideals and
$I_k=J_k$ for $k \leq 5$ and $I_k = J'_k$ for $k \geq 5$.
Hence $I$ has a segmentwise critical Hilbert function.
However, for any linear transformation $\varphi$ on $S$,
if $\phi(I)_4$ is a canonical critical space,
then $\phi(I)_6$ is not a canonical critical space.
}
\end{Example}

Finally, we generalize Peeva's result \cite[Corollary 1.4]{Peeva}.
For the $(n-1)$-th Macaulay representation
$a= \sum_{j=1}^s { a(j) + n-j  \choose n-j}$ of a positive integer $a>0$,
define
\begin{eqnarray*}
a^+ =
\left\{
\begin{array}{ll}
a \Mac, &\mbox{ if } s \geq n-2,\\
a \Mac + {n-s-1 \choose n-s-1} + \cdots + {2 \choose 2}, &
\mbox{ if } s < n-2.
\end{array}
\right.
\end{eqnarray*}

\begin{Theorem}
\label{Peeva}
Let $I$ be a homogeneous ideal
satisfying, for all $k \geq 1$,
$$
\mbox{if $H(I,k) \Mac < H(I,k+1)$ then $H(I,k+2) \leq H(I,k+1)^+$}.
$$
Then $I$ is Gotzmann.
\end{Theorem}

\begin{proof}
We show that $I_k$ is Gotzmann for all $k$
by using induction on $k$.
It is clear that $I_k$ is Gotzmann for $k \gg 0$.
We show if $I_{k+1}$ is Gotzmann then $I_k$ is Gotzmann.

Let
$H(I,k)= \sum_{j=1}^s { k-a_j + n-j \choose n-j}$
be the $(n-1)$-th Macaulay representation of $H(I,k)$.
If $H(I,k) \Mac = H(I,k+1)$
or
$H(I,k-1) \Mac = H(I,k)$ then $I_k$ is Gotzmann by Lemma \ref{obvious}.
Thus we may assume $H(I,k) \Mac < H(I,k+1)$ and $H(I,k-1) \Mac < H(I,k)$.
Then, by the assumption, we have
$s < n-2$ and $H(I,k+1) \leq H(I,k)^+$.
Hence
$$H(I,k+1)= \sum_{j=1}^s { k+1 -a_j +n-j \choose n-j} + {n-s-1 \choose n-s-1} + \cdots + {n-\ell \choose n-\ell}$$
for some $\ell \leq n- 2$.
Since $I_{k+1}$ is Gotzmann by the induction hypothesis,
Gotzmann's persistence theorem implies
$$H(I_{\langle k+1 \rangle},t)
= \sum_{j=1}^s { t -a_j +n-j \choose n-j} 
+ \sum_{j=s+1}^{\ell} {t-(k+1)+n-j \choose n-j}
\ \ \mbox{ for all } t \geq k+1.$$
Since $n-\ell \geq 2$,
by Proposition \ref{Gotzmann1989},
which we will show later in the appendix,
it follows that $(I_{\langle k+1 \rangle})^\sat$ is a critical ideal of type
$(a_1,\dots,a_s,k+1,\dots,k+1)$.
Then, since 
$H(I,k)=\sum_{j=1}^s { k-a_j + n-j \choose n-j}
=H((I_{\langle k+1 \rangle})^\sat,k)$
and $I_k \subset (I_{\langle k+1\rangle})^\sat$,
we have $I_k = ((I_{\langle k+1 \rangle})^\sat)_k$.
Hence $I_k$ is a critical space, and therefore $I_k$ is Gotzmann as desired.
\end{proof}

\begin{Example}
{\em
By Macaulay's theorem, there is a one-to-one correspondence between Hilbert functions of
homogeneous ideals of $A$ and lexsegment ideals of $A$.
Let 
$I=(u_1,u_2,\dots, u_t)$
be a lexsegment ideal.
Set $d_k =\deg u_k$ for $k=1,2,\dots,t$.
Suppose $d_1 \leq d_2 \leq \cdots \leq d_t$.
If $A=K[x_1,x_2]$,
then segmentwise critical Hilbert function corresponds to lexsegment ideal $I$
satisfying $d_2  < d_3 -1 < d_4 -2 < \cdots <d_t -t+2$
and Hilbert functions satisfying the condition of Theorem \ref{Peeva} correspond to
lexsegment ideals $I$ satisfying $d_1=\cdots =d_{s_1} < d_{s_1+1}-1 = \cdots =d_{s_2}-1 < d_{s_2+1} -2 = \cdots$,
where $1 \leq s_1 < s_2 < \cdots$.
Note that if $n=2$ then Theorem \ref{Peeva} is equivalent to Peeva's result
since $a^+= a \Mac$.
}\end{Example}

\begin{Remark}
{\em
We can construct new inflexible Hilbert functions
by using Theorems \ref{2-5} and \ref{Peeva}.
Let $d>0$ be an integer and
$H$ a segmentwise critical function such that $H(k) \Mac =H(k+1)$ for $k \geq d$.
Let $H'$ be a Hilbert function of a homogeneous ideal of $A$ satisfying the condition of Theorem \ref{Peeva}.
If $H(d)=H'(d)$ and $H(d+1)=H'(d+1)$,
then we define the new Hilbert function $\tilde H$ by
$\tilde H(k)=H(k)$ for $k \leq d$
and $\tilde H(k)=H'(k)$ for $k \geq d$.
This new function $\tilde H$ is an inflexible Hilbert function.

Indeed, if a homogeneous ideal $I$ of $A$ has the Hilbert function $\tilde H$,
then, by the same way as in the proof of Theorem \ref{Peeva},
it follows that $I_k$ is a Gotzmann space for $k \geq d$.
On the other hand, by Gotzmann's persistence theorem,
the ideal $J= \sum_{k=0}^{d} I_{\langle k \rangle}$ has the Hilbert function $H$.
Since $H$ is a segmentwise critical function,
$I_k=J_k$ is a Gotzmann space
for $k \leq d$.
}
\end{Remark}

\begin{Example}
{\em
Let $d >0$ be an integer and
let $I=(x_1^d) + (x_1,\dots,x_n)^{d+1} \subset A$.
Then all ideals of $A$ having the same Hilbert function as $I$
are Gotzmann.
However, the Hilbert function of $I$ is not segmentwise critical
and does not satisfy the condition of Theorem \ref{Peeva}.
}
\end{Example}

\section*{Appendix}
\renewcommand{\thesection}{\Alph{section}}
\setcounter{section}{1}
\setcounter{Theorem}{0}

Let $I$ be a homogeneous ideal of
$A = K[x_1, \ldots, x_n]$ and
write $P_I(t)$ for the Hilbert polynomial of $I$.
It follows from (\ref{MacEq})
that $P_I(t)$ can be written uniquely
in the form
\begin{eqnarray}
\label{hilbertpolynomial}
P_I(t) = {t - a_1 + n - 1 \choose n - 1}
+ \cdots +
{t - a_s + n - s \choose n - s},
\end{eqnarray}
where $1 \leq s \leq n - 1$ and where $a_1, \ldots, a_s$
are integers with
$0 < a_1 \leq \cdots \leq a_s$.

Gotzmann \cite[Proposition 2]{Got}
proved using the languages of algebraic geometry
that if a homogeneous ideal $I$ of $A$ is saturated
and if the Hilbert polynomial of $I$
of the form $(\ref{hilbertpolynomial})$
satisfies $n - s \geq 2$, then
there exists a linear transformation
$\varphi$ on $A$
such that $\varphi(I)$ is a canonical
critical ideal.
Gotzmann \cite{Got1982,Got} also show that
the Hilbert scheme parametrizing
ideals with Hilbert polynomial (\ref{hilbertpolynomial}),
which satisfies either 
(i) $n-s \leq 2$ or (ii) $s=n-1$ and $a_1- a_n \leq 1$,
is irreducible and simply connected, with a single Borel-fixed point.
On the other hand,
Mall \cite{Mall} gave a combinatorial proof of the simple connectivity
of these Hilbert schemes
by showing that they have a single Borel-fixed point.

Theorem \ref{gotzmann} together with
Proposition \ref{Gotzmann1989}
stated below yields a
simple and purely algebraic proof of \cite[Proposition 2]{Got}.

\begin{Proposition}[Gotzmann]
\label{Gotzmann1989}
Let $I$ be a saturated ideal of $A$
with the Hilbert polynomial $(\ref{hilbertpolynomial})$
and suppose that one of the following conditions
is satisfied:
\begin{enumerate}
\item[(i)]
$n - s \geq 2$;
\item[(ii)]
$s = n - 1$ and
$a_{n-1} - a_1 \leq 1$.
\end{enumerate}
Then the Hilbert function of $I$ is the critical function of type
$(a_1, a_2, \ldots, a_s)$.
\end{Proposition}

Note that a universal lexsegment ideal whose Hilbert function
is a critical function of type $(a_1,\dots,a_s)$
with $n- s \geq 2$ has at most $n-2$ generators.
In particular, by (\ref{unilex}) such a universal lexsegment
ideal has no generators which are divisible by $x_{n-1}$.

Proposition \ref{Gotzmann1989} was proved in \cite{Got1982} but
we include a proof for the sake of completeness.
For the $d$-th Macaulay representation
$a= \sum_{j=k}^d { a(j) + j\choose j}$
of a positive integer $a>0$,
let
\[
a_{\langle d \rangle}
= {a(d) + d - 1 \choose d}
+ \cdots +
{a(k) + k - 1 \choose k}
\]
and $0_{\langle d \rangle}=0$.
Then, by (\ref{MacIneqSpace}) and (\ref{MacEq}),
any homogeneous ideal $I$ of $S$ satisfies
\begin{eqnarray}
\label{fundamental}
H(I,k+1)_{\langle n-1 \rangle } \geq H(I,k),
\, \, \, \, \, \, \, \, \, \,
k = 0, 1, 2, \ldots.
\end{eqnarray}

\begin{proof}[Proof of Proposition \ref{Gotzmann1989}]
We define functions
$H : \ZZ_{\geq 0} \to \ZZ_{\geq 0}$ 
and $H' : \ZZ_{\geq 0} \to \ZZ_{\geq 0}$
by setting
\[
H(t) = {t - a_1 + n - 1 \choose n - 1}
+ \cdots +
{t - a_s + n - s \choose n - s}
\]
and
\begin{eqnarray}
\label{a0}
&&H'(t) = {t - a_1 + n - 2 \choose n - 2}
+ \cdots +
{t - a_s + n - s -1\choose n - s-1}.
\end{eqnarray}
Then a simple computation yields
\begin{eqnarray}
\label{a1}
H'(t)=H(t) -H(t-1)
\ \ \mbox{ for all $t \geq 0$}.
\end{eqnarray}

By considering an extension field,
we may assume that $K$ is infinite.
Then, since $I : (x_1, \ldots, x_n) = I$,
it follows that
the quotient algebra $A / I$
possesses a homogeneous
non-zero-divisor
$\theta$ with $\deg \theta = 1$.
Let $B = K[y_1, \ldots, y_{n-1}]$
denote the polynomial ring
in $n - 1$ variables over $K$.
Regarding
$(A / I)/\theta(A / I)$
as a quotient algebra
of $B$ by a homogeneous ideal
$J$, a simple computation yields
that
\[
H(B/J, t) = H(A/I,t) - H(A/I,t-1)
\ \ \ \  \mbox{ for all $t \geq 0$}.
\]
Then, since $ H(B,t) + H(A,t-1) = H(A,t)$,
we have
\begin{eqnarray}
\label{a3}
&& H(J, d) = H(I,t) - H(I,t-1)
\ \ \mbox{ for all $t \geq 0$}.
\end{eqnarray}

Since the Hilbert polynomial
$P_I(t)$ of $I$ is $(\ref{hilbertpolynomial})$,
it follows that $H(I,t) = P_I(t) = H(t)$
for all $t \gg 0$.
Then we have $H(J,t)=H'(t)$ for all $t \gg 0$ by
(\ref{a1}) and (\ref{a3}).
We claim that if either (i) or (ii) is satisfied then
\begin{eqnarray}
\label{aX}
&& H(J,t) \leq H'(t)
\ \ \mbox{ for all $t \geq 0$}.
\end{eqnarray}

\textit{Case 1}:
Suppose $n-s \geq 2$.
Then (\ref{a0}) coincides with the $(n-2)$-th Macaulay representation
of $H'(t)$.
Thus $H'(t+1)_{\langle n-2 \rangle} = H'(t)$
for all $t \geq 0$.
Now, if $H(J,t_0+1) \leq H'(t_0+1)$ for some $t_0$,
then it follows from $(\ref{fundamental})$
that
\[
H(J,t_0)
\leq H(J,t_0+1)_{\langle n-2 \rangle}
\leq H'(t_0+1)_{\langle n-2 \rangle}
= H'(t_0).
\]
(The second inequality follows from Lemma \ref{macaulayexpand}.)
Hence
$H(J,t_0) \leq H'(t_0)$.
Consequently,
we have
$H(J,t) \leq H'(t)$
for all $t$ as desired.

\textit{Case 2}:
Suppose that
$s = n - 1$,
$a_{n-1} - a_1 \leq 1$
and
$a_p = a_{p+1} = \cdots = a_{n-1} $,
where either $a_{p-1} < a_p$ or $p = 1$.
Then it follows from Lemma \ref{2-0} that
\begin{eqnarray*}
&&H'(t) = \sum_{k=1}^{p-1} {t -a_k + n-1-k \choose n-1-k}
+{t-(a_p-1) +n-1-p \choose n-1-p}
\end{eqnarray*}
is the $(n-2)$-th Macaulay representation of $H'(t)$
for all $t \geq a_p$.
Thus $H'(t+1)_{\langle n-2 \rangle} = H'(t)$ for all $t \geq a_p$.
Then, by the same way as in Case 1,
we have $H(J,t) \leq H'(t)$
for all $t\geq a_p-1$.
On the other hand, since $H(J,a_p-1) \leq H'(a_p-1)< n-1$,
it follows that $H(J,t)=0$ for $t < a_p-1$.
Thus the desired inequality
(\ref{aX}) follows.
\medskip

Now (\ref{a1}), (\ref{a3}) and (\ref{aX})
imply
\begin{eqnarray}
\label{hoshihoshi}
&&H(I,t)= \sum_{k=1}^t H(J,t) \leq \sum_{k=1}^t H'(t) =H(t)
\ \ \mbox{ for all $t \geq 0$}.
\end{eqnarray}
Since $H(I,t)=P_I(t)=H(t)$ for all $t\gg 0$,
(\ref{hoshihoshi})
implies $H(J,t)=H'(t)$ for all $t \geq 0$.
Hence, by (\ref{hoshihoshi}),
we have $H(I,t)=H(t)$ for all $t \geq 0$.
Thus the Hilbert function of $I$
is the critical function of type $(a_1,\dots,a_p)$ as desired.
\end{proof}

\begin{Remark}
{\em
If $I$ has a critical Hilbert function of type $(a_1,\dots,a_p)$ with
$p \leq n-1$, then the Hilbert polynomial of $I$ is given by the formula
(\ref{hilbertpolynomial}).
However, this does not imply $H(I,k)=P_I(k)$ for all $k \geq 0$.
Indeed, if $n=3$, $p=1$ and $a_1=4$ then
$P_I(t)= { t-2 \choose 2}= \frac{(t-2)(t-3)} 2$.
Thus $P_I(1)=1$, however, $H(I,1)= { -1 \choose 2}=0$.
}
\end{Remark}

\begin{Remark}
{\em
Proposition \ref{Gotzmann1989}
is false if neither (i) nor (ii) is satisfied.
Indeed, it follows from
\cite[Theorem 26]{Mall} that
if a sequence of integers
$0<a_1 \leq \cdots \leq a_s$ with $1 \leq s \leq n-1$
satisfies neither (i) nor (ii), then there exists a
saturated non-lexsegment Borel-fixed ideal with the Hilbert polynomial 
(\ref{hilbertpolynomial}).
The Hilbert function of such an ideal is not critical by
Lemma \ref{gin}.
}
\end{Remark}

\noindent
\textbf{Acknowledgments}:
We would like to thank the referee
for suggesting the study of segmentwise critical functions
as well as for many helpful comments to improve the paper.

\end{document}